\documentclass[12pt,reqno,twoside]{amsart}

\usepackage{amscd}
\usepackage{amsmath}
\usepackage{amssymb}
\usepackage{amsthm}
\usepackage{epsf}
\usepackage{verbatim}
\usepackage{graphicx}
\usepackage{parskip}
\usepackage{caption,hyperref}
\usepackage{color,comment}
\input epsf.tex

\usepackage{etoolbox}
\makeatletter
\patchcmd{\@maketitle}
  {\ifx\@empty\@dedicatory}
  {\ifx\@empty\@date \else {\vskip0ex \centering\footnotesize\@date\par\vskip1ex}\fi
   \ifx\@empty\@dedicatory}
  {}{}
\patchcmd{\@adminfootnotes}
  {\ifx\@empty\@date\else \@footnotetext{\@setdate}\fi}
  {}{}{}
\makeatother

\theoremstyle{theorem} 

\newtheorem*{theorem*}{Theorem}

\theoremstyle{definition} 

\theoremstyle{definition}

\theoremstyle{definition}

\theoremstyle{definition}

\theoremstyle{definition}

\usepackage{chngcntr}
\counterwithout{figure}{section} 

\begin{document}

\title[Assessing congressional districting  in ME and NH]{Assessing congressional districting\\ in Maine and New Hampshire}

\author[]{Sara Asgari, Quinn Basewitz, Ethan Bergmann, Jackson Brogsol, \\ Nathaniel Cox, Diana Davis, Martina Kampel, Becca Keating,\\  Katie Knox, Angus Lam, Jorge Lopez-Nava, Jennifer Paige, \\Nathan Pitock,  Victoria Song, Dylan Torrance}

\date{\today}

\begin{abstract}
We use voting precinct and election data to analyze the political geography of New Hampshire and Maine. We find that the location of dividing line between Congressional districts in both states are significantly different than what we would expect, which we argue is likely due to incumbent gerrymandering. We also discuss the limitations of classical fairness measures for plans with only two districts.
\end{abstract}

\maketitle

\section{Our work}

\begin{figure}[!h]
\begin{center}
\includegraphics[height=0.4\textwidth]{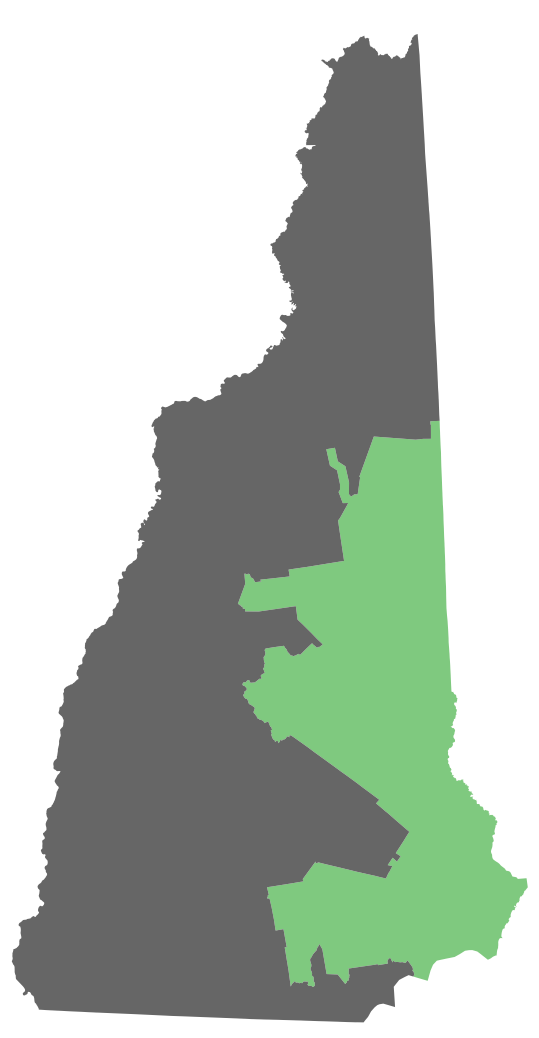} \hspace{0.5in}
\includegraphics[height=0.4\textwidth]{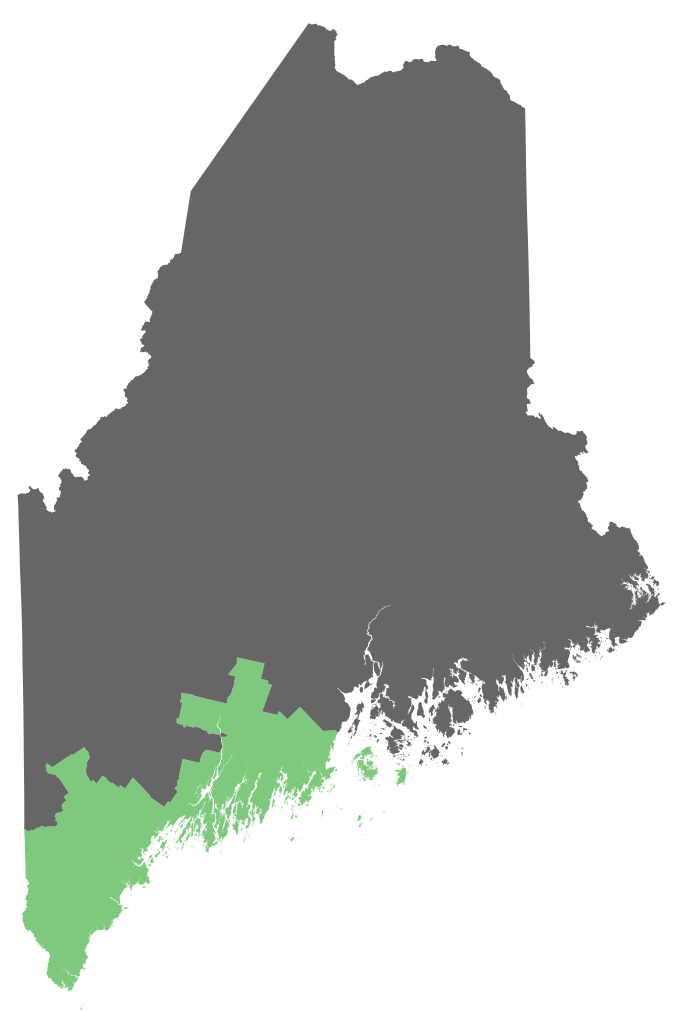} 
\caption{The congressional districts of New Hampshire and Maine enacted in 2011. They have equal population according to the 2010 census. \label{fig:current_districts}}
\end{center}
\end{figure}

A guiding principle of representative democracy is \emph{one person, one vote}, the philosophy that everyone's vote counts equally. The practice of political \emph{gerrymandering} seeks to dilute the voting power of groups of people, by making it more difficult for them to convert their votes into representative seats. We assessed the fairness of the 2011 congressional districting maps in New Hampshire and Maine (Figure \ref{fig:current_districts}), by measuring the enacted plan against a background distribution of thousands of other legal plans.

We downloaded voting precinct maps from each state, and election results reported by precinct. We then created thousands of legal districting plans for each of the two states: connected, and equal in population. Finally, we assessed the properties of the enacted plan against the background distribution of the associated properties of the comparison plans.

\begin{figure}[!h]
\begin{center}
\includegraphics[height=0.5\textwidth]{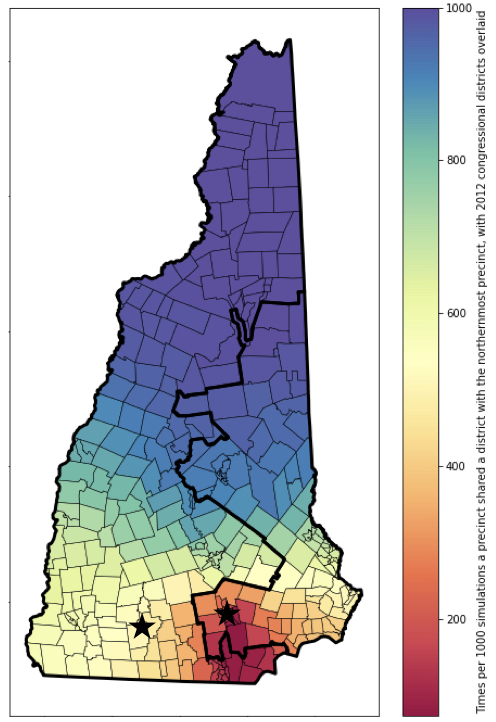}  \hspace{0.5in} \includegraphics[height=0.5\textwidth]{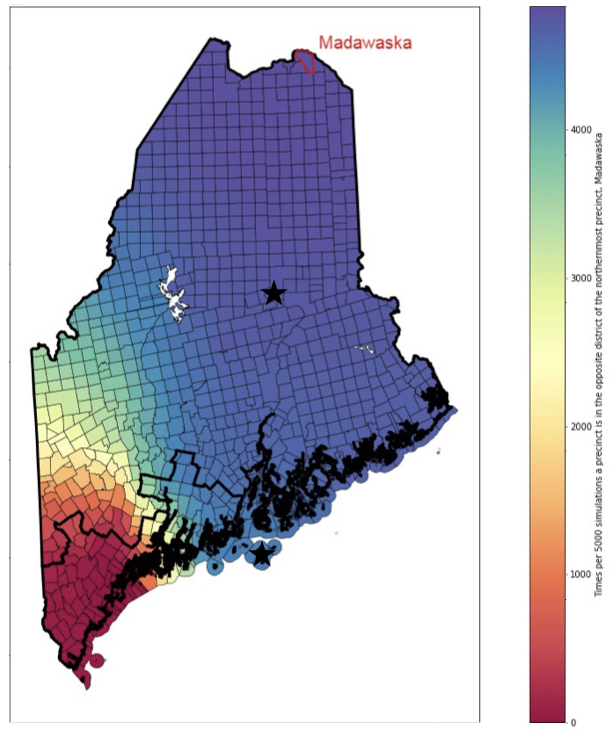}
\caption{The enacted congressional districts lines (black) compared to the typical dividing areas (yellow). Residence locations of the incumbent representatives are shown as black stars. \label{fig:rainbow_districts}}
\end{center}
\end{figure}

One striking finding is shown in Figure \ref{fig:rainbow_districts}. The yellow region shows where equal-population districting plans tend to divide the state. In both New Hampshire and Maine, the division between the enacted districts (black) is roughly perpendicular to where we would expect it to be (yellow). Why is this?

A plausible explanation is \emph{incumbent gerrymandering}. The black stars show where the incumbent representatives lived in 2011, when the current districting plan was designed. With the enacted division, they are in different districts; with a division in the yellow region, they would be in the same district, and thus competing against each other in the following election, an undesired event. We found this result in both New Hampshire and Maine.

\section{Background on districting and gerrymandering}

As dictated in the U.S. Constitution, a national census is taken every ten years, most recently in 2020, to count and record the geographic distribution of the population. The census data are used to re-allocate federal funding, inform other community work, and importantly, re-allocate congressional representation. The Constitution mandates that districts have equal population based on the most recent census, so after each census the state's Congressional districts are redrawn. The new maps are generally proposed by a legislative board or an outside committee, and passed as a bill by the state house of representatives. 
The current congressional districts for Maine and New Hampshire, enacted in 2011 following the 2010 census, are shown in Figure \ref{fig:current_districts}.

Politicians often use a technique known as \emph{gerrymandering} in order to manipulate boundaries within state or congressional maps for purposes of partisan gain. The party holding the majority of seats at the state level has the power to approve the maps, and sometimes draws districts in a way that dilutes the voting power of the opposing party. In the past, the \emph{compactness} of a district, measured by e.g. the ratio of its border length to its area, has been used as a measure of the degree to which a district is gerrymandered. For example, one might observe that the boundaries between districts in Figure \ref{fig:current_districts} wander back and forth and have many protrusions, and posit that this is evidence of nefarious intentions. But with modern computing and data analysis, even plump, reasonable-looking districts can be customized to achieve political outcomes. We will describe techniques using random walks and outlier analysis that detect gerrymandering based on the political \emph{effects} of the districting plan, rather than its geometric \emph{shape}. 

Gerrymandering techniques, like \emph{packing} opponents into a district that they win by a landslide and \emph{cracking} the rest into districts that they just barely lose, are commonly used by the incumbent party, and can result in a governing body that significantly misrepresents the constituents. Politicians use gerrymandering tactics to increase their chances of reelection, muting those who would provide opposition. Each state has different requirements for its districts, and conducts its own re-mapping procedure independently following the U.S. census every 10 years, and there are no detailed federal guidelines for redistricting, so there are currently few regulations in place to prevent incumbent politicians from gerrymandering.

\section{Political background}

\subsection{Maine}

Maine is divided into just two congressional districts, and is thus represented by two members in the U.S House of Representatives. This has not always been the case: at its peak, Maine was represented by eight members in the U.S. House; as late as the 1950s, Maine had three congressional districts. It was only beginning in 1962 that the current two-state district allocation came into effect \cite{Maine-Encyclopedia}.

Maine’s 1st Congressional District (ME-1) is the geographically smaller district, covering the southern coastal areas – and more populous towns, cities, and counties – of the state. It includes the city of Portland, as well as major towns like Brunswick. The 2nd Congressional District (ME-2), on the other hand, comprises nearly 80 percent of the state’s total land area, but is sparsely populated. It includes Aroostook County, which is almost two-thirds the size of New Hampshire and larger than the state of Connecticut but has only 68,000 people \cite{Maine-Encyclopedia}.

Over the past several decades, ME-2 has continued to expand geographically – likely a response to the state’s population trends shifting toward the southern coastal cities. As mandated by the U.S. Constitution, congressional districts must be equal in population based on the most recent census. But between 2000 and 2009, for example, Aroostook County (ME-2) saw a population decline of 3.3\%, while Cumberland County (ME-1) saw a population increase of 4.9\% \cite{Census}, so by 2009, ME-2 voters had more voting power than those in ME-1. With economic growth becoming more concentrated in the state’s urban areas (in ME-1), and with the state’s economy shifting away from manufacturing and toward the healthcare and finance industries \cite{Colgan}, these population trends will likely continue in the coming decades.

The economic disparities between ME-1 and ME-2, hinted above, have led to the idea of “Two Maines” – that is, a perception that while the coastal and southern counties of ME-1 are relatively prosperous, the rural north and western counties are economically depressed or backwards \cite{Spruce}. Examining the census data shows that there is some truth to this: the median household income in ME-1 is \$64,599, compared to \$48,603 in ME-2; the poverty rate in ME-1 is 8.7\%, compared to 13.6\% in ME-2. In terms of educational attainment: 38.2\% of residents of ME-1 have earned a bachelor’s degree or higher, while only 25.6\% of the residents of ME-2 have done the same \cite{acs}.

The politics of the two congressional districts has reflected, to some extent, the economic and cultural differences between northern and southern Maine. ME-1, with its urban areas, has had more of a liberal streak when it comes to federal, statewide, and local elections. The more rural ME-2, on the other hand, has more often leaned conservative. In the 2016 presidential race, for example, Republican candidate Donald Trump managed to beat Democratic candidate Hillary Clinton in  ME-2 (winning by 11 points), even as Clinton won a majority and plurality of votes in ME-1 and the statewide vote, respectively \cite{NYT}. The Cook Political Report, in its most recent 2017 report, gave ME-1 a partisan voting index score of D+8, while giving ME-2 a score of R+2 \cite{Cook}.

Then again, Mainers often describe themselves as free from partisan influences and fiercely independent – a claim bolstered by Maine’s history of electing independent political officials, from then-Governor James Longley, to then-Governor and now-U.S. Senator, Angus King. Indeed, according to recent surveys, “unaffiliated,” “unenrolled,” and “independent” voters almost outnumber the number of Democrats in the state, and greatly outnumber Republicans. Out of the 1,063,383 registered voters in Maine (as of July 14, 2020), 339,782 registered as “unenrolled” and 41,693 registered as “green independents,” while 386,786 registered as Democrats and 295,122 registered as Republicans \cite{Voters}. In short, of all the registered voters in Maine, around 36\% remain unaffiliated from the two major parties – around the same percentage as those who registered as Democrats, and about 8 percentage points more than those who registered as Republicans.


\subsection{New Hampshire}\label{nh-intro}

New Hampshire is the 10th least populous of the 50 US states, following Maine. When the United States Congress was founded in 1789, New Hampshire had three representatives, but it is currently divided into just two congressional districts. As a swing state, and as the first state to vote in the Presidential primary, New Hampshire is strongly targeted by campaigns in competitive elections. 

New Hampshire’s 1st Congressional District (NH-1) is the smaller district geographically and includes Manchester -- New Hampshire's largest city -- and the seacoast and the lakes region. The 2nd Congressional District (NH-2) is larger, and includes areas that are thinly populated, such as the northernmost Coos county adjacent to Canada, while also including New Hampshire's second- and third-largest cities, Nashua and Concord.

The first district has a population of 673,194, with the second district slightly smaller at 669,601. Unlike the disparities seen between Maine's two districts, New Hampshire's two districts are much more economically similar. The median household income in NH-1 is \$73,488, comparable to \$73,249 in NH-2. One difference worth noting is the proportion of residents in rural versus urban locations. In NH-1, 70\% of the residents live in urban areas; in contrast, in NH-2, only 51\% are located urban areas.

In New Hampshire, there is a long-standing tradition of registering as an independent rather than with the Democratic or Republican Party. Although voters registered as independent cannot vote in party primaries, it is common for voters to change their registration to participate in the primary election and change their status back to independent after casting their ballot. In addition, New Hampshire is the home to the Free State Project, which recruits libertarians to move to New Hampshire in order to make the state a stronghold for libertarian ideas. All told, 42\% of New Hampshire voters are not registered with the Democratic or Republican Parties. In general elections, many of these independent voters choose to cast their ballot for a Democrat or a Republican, while a substantial number vote for a third-party candidate. Independent voters have the power to decide elections and swing results in favor or one party over another.

\section{Gerrymandering background}\label{gerry-background}

To analyze the fairness of an enacted districting plan, we generate a large ensemble of valid districting plans, and then compare the enacted plan against this background distribution, a method known as \emph{outlier analysis} \cite{duchin}. Representing the population units of a state (e.g. census blocks, voting precincts) as graph vertices, we connect two vertices with an edge if they represent neighboring units. We then partition the state into districts, by grouping vertices of the graph that comprise connected subgraphs that have almost equal populations and follow the state guidelines. These guidelines might limit splitting the counties and cities, or having unusually-shaped or noncontiguous districts. 

To divide the graph into districts that follow the above-mentioned rules, we used the GerryChain software library with the recombination (``ReCom'') method \cite{GerryChain,ReCom}. The ReCom method works as follows: starting with a districting plan, it randomly selects two adjacent districts and merges them. It then finds a random spanning tree\footnote{\emph{Spanning tree}: a connected acyclic graph containing all of the graph's vertices.} of the merged subgraph using Dijkstra's algorithm\footnote{ Dijkstra's algorithm is designed to return a \emph{minimum-weight spanning tree}. In order to use the algorithm to find a random spanning tree, GerryChain assigns each graph edge a random number between 0 and 1, and the algorithm returns the spanning tree with the minimum sum of all edge weights.} and splits the tree into two trees of the same size by cutting a graph edge. The subgraphs represented by these trees comprise our new districts. If it is not possible to cut the tree into two equal-size districts, or if the new districts do not follow the state's guidelines, the method starts over by finding another random spanning tree. The method continues the process until it creates the required number (typically thousands or millions) of districting plans.

\begin{figure}[h]
\centerline{\includegraphics[scale=0.4]{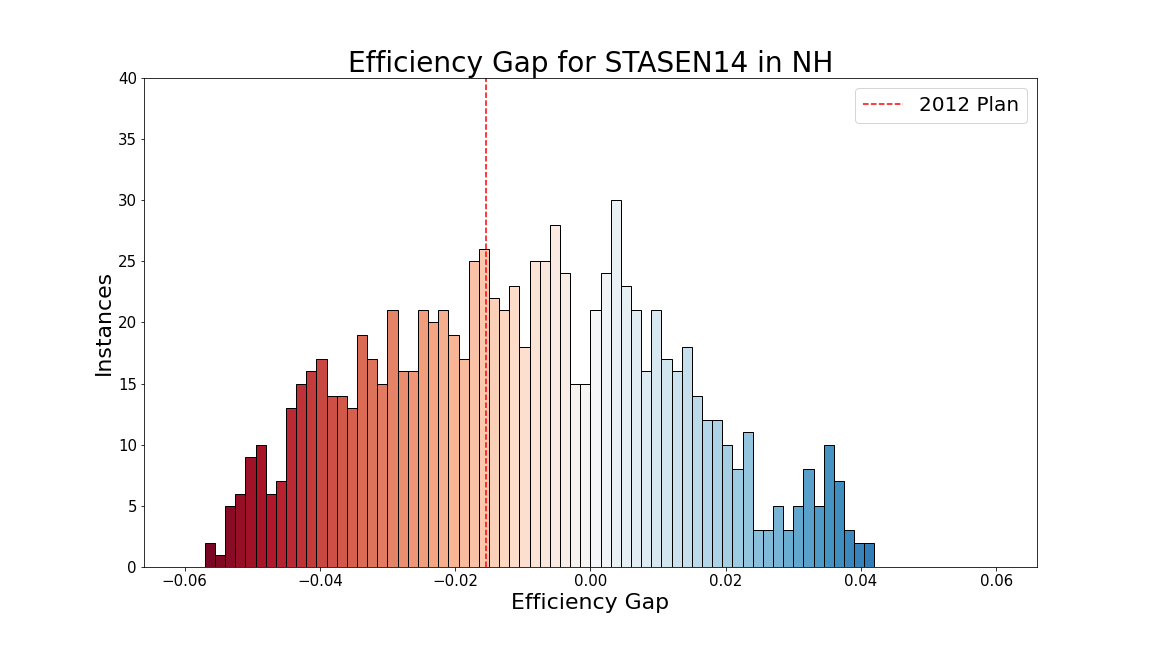}}
\caption{Histogram of the efficiency gap among an ensemble of 1000 plans for the 2014 State Senate Election in New Hampshire. The score for the 2012 plan is marked with a vertical dashed line.}
\label{NHEG}
\end{figure}

Once we have an ensemble of districting plans, we can analyze whether the enacted plan is an outlier compared to the others. We first measure a ``fairness'' score (e.g. the mean-median, see \S\ref{sec:measures}) for each plan in the ensemble. 
 We create a histogram, with the fairness scores along the horizontal axis and the number of plans that achieve a given range of scores along the vertical axis. As an example, the graph in Figure \ref{NHEG} uses the efficiency gap, and the enacted plan is not an outlier since it falls in the middle of the scores for the ensemble. 

\section{Measures of fairness in a state with two districts}\label{sec:measures}
Substantial literature exists outlining methods to evaluate partisan gerrymandering in a variety of situations. The primary mechanisms for assessing gerrymandering are mean-median distribution, partisan bias, partisan symmetry, efficiency gap, and sampling \cite{duchin}. These methods, termed ``measures of fairness,'' are informative but not suitable for all cases. Specifically, these mechanisms become less informative as the number of districts decreases. Measures of fairness in states with two districts (namely, New Hampshire and Maine) are particularly uninformative.

The authors of \cite{paradoxes} explain two common measures of fairness:

\begin{itemize}
    \item[] The \emph{mean-median} metric is vote-denominated: it produces a signed number that is often described as measuring how far short of half of the votes a party can fall while still securing half the seats. A similar metric, \emph{partisan bias}, is seat-denominated. Given the same input, it is said to measure how much more than half of the seats will be secured with half of the votes. The ideal value of both of these metrics is zero.
\end{itemize}

\begin{figure}[h]
\centerline{\includegraphics[scale=0.4]{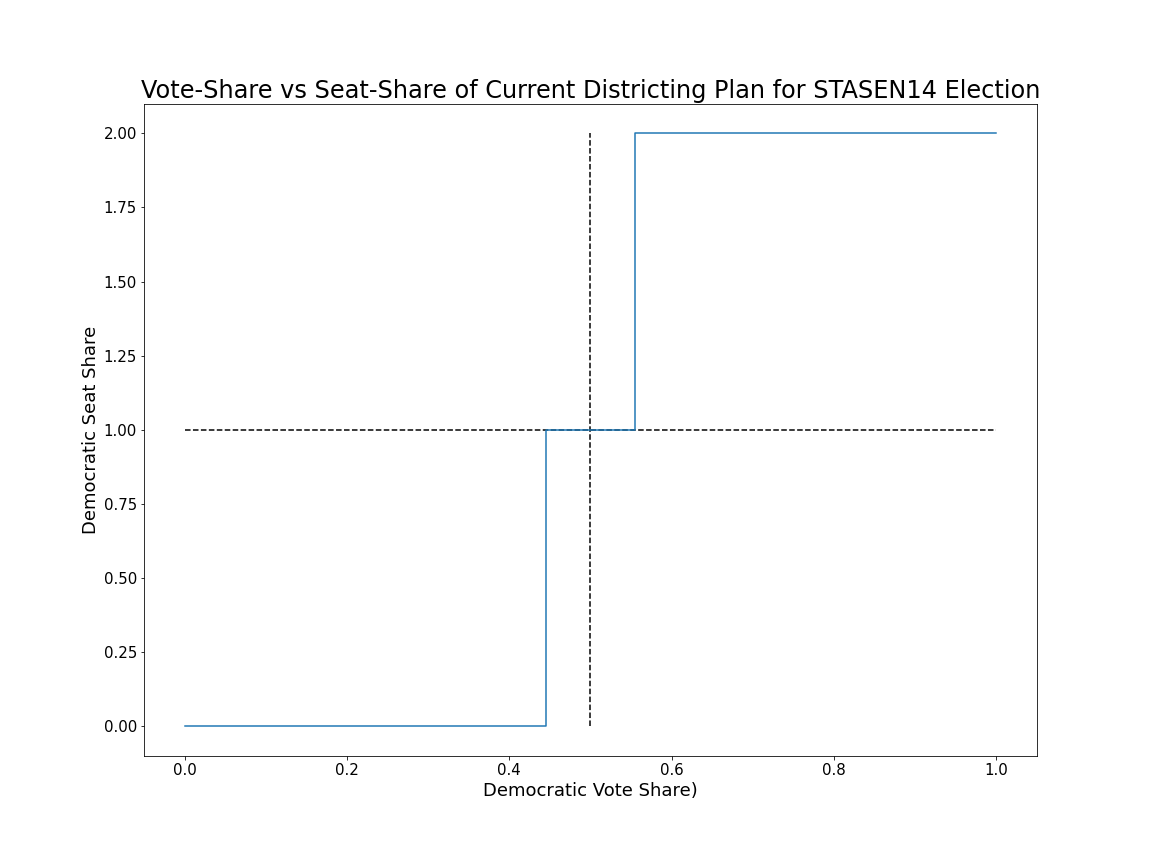}}
\caption{Vote shares versus seat shares plot for the 2014 State Senate Election in New Hampshire}
\label{NHVSSS}
\end{figure}

In a state with two districts, the value of the mean-median metric and partisan bias metric is always zero (and therefore unenlightening). To explain this phenomenon, Figure \ref{NHVSSS} shows the Vote-Share vs Seat-Share curve (blue) of the New Hampshire 2014 State Senate Election. 

The mean-median score is the horizontal distance from the blue curve to the middle point (0.5, 1.0), and the partisan bias score is the vertical distance from the curve to the same middle point. Regardless of where the curve flips in a two-district plan, the curve passes through the central point, so both the mean-median and partisan bias scores are 0, making it an uninformative metric for analysis in New Hampshire, Maine and many other states.

Another metric, \emph{efficiency gap}, can be somewhat enlightening but has flaws as well. Efficiency gap is defined as
\begin{equation*}
   \frac{\text{wasted Rebuplican votes} - \text{wasted Democratic votes}}{\text{total votes}}.
\end{equation*}
The optimal efficiency gap is zero, as this implies that voters in both parties waste the same proportion of their votes. 

\begin{figure}[h]
\centerline{\includegraphics[scale=0.4]{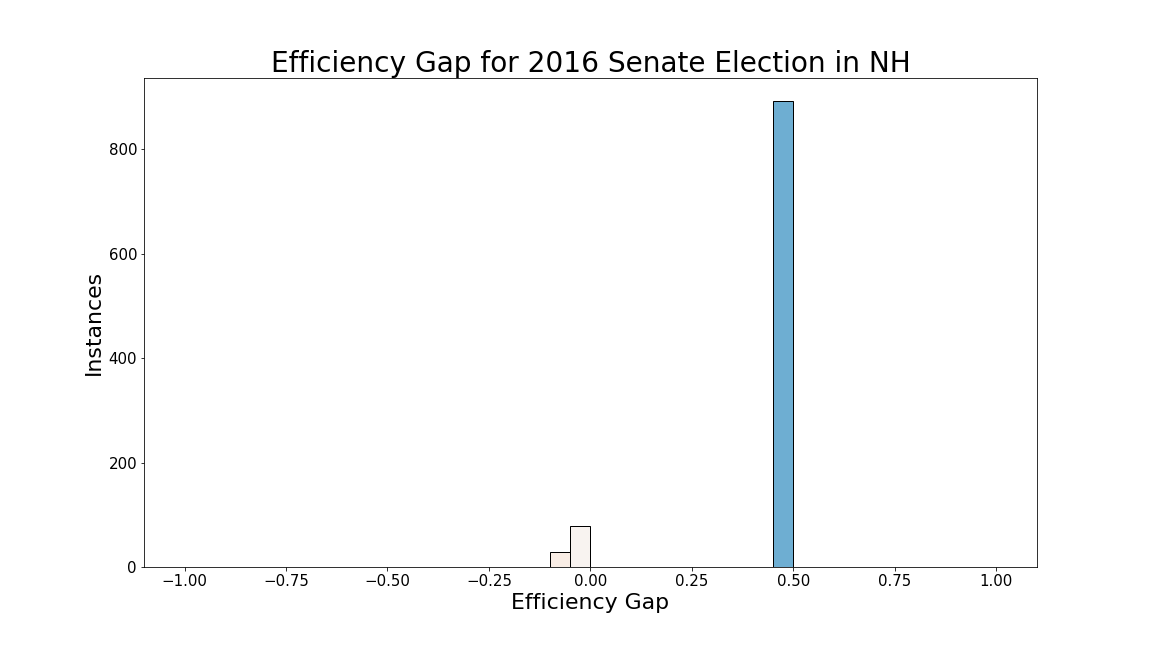}}
\caption{Histogram of the efficiency gap among an ensemble of 1000 plans for the 2016 Senate Election in New Hampshire.}
\label{NHBadEG}
\end{figure}

In a two-district election, the efficiency gap is always close to 0 or 0.5, depending on election dynamics; an example of this phenomenon is in Figure \ref{NHBadEG}. For a given election, for all of the districting plans where one party wins both seats, the efficiency gap is identical for every plan, making it useless for creating a nice background distribution histogram like the one in Figure \ref{NHEG}. If each party wins one of the two seats in most or ideally all of the districting plans in the ensemble, then the efficiency gap histogram is useful for our analysis. 

As a result of the uninformative nature of mean-median and partisan fairness scores in a two-district state, for New Hampshire and Maine, our analysis used the efficiency gap. Our results are meaningful, but limited because the efficiency gap struggles to provide strong evidence in states with few districts. Further work could develop new measures of fairness in cases where traditional measures fail to be meaningful.

\section{Data gathering and cleaning}

\subsection{New Hampshire}
New Hampshire election data is available on the New Hampshire Secretary of State website. In New Hampshire, each town is a voting precinct, and cities are broken up into voting precincts called \emph{wards}. The statewide precinct map shape file is available on NH GRANIT: New Hampshire's Statewide Geographic Information System Clearinghouse. 

We had to edit the data slightly in order for the precincts in the shape file to match exactly with the election data. The shape of the precinct \emph{Somersworth - Ward 3} is made of two polygons that meet at a corner, which the shapefile included as two pieces; we merged them. For the precinct \emph{Derry}, election data was reported for the entire town, while the shapefile included five separate wards. We merged the wards in order for the election data and shapefile to match.

\subsection{Maine}
Like New Hampshire, Maine's election data is available on the Secretary of State website, and the precinct map shape file was available on Maine’s geolibrary site. 

Voting data in Maine is most frequently recorded by town, as (like New Hampshire) nearly every town is its own precinct. But while our precinct shapefile contained every precinct in Maine, each election's results only reported precincts with votes. This would not be a problem in more populated states, but Maine is largely uninhabited:  over one-third of the 924 precincts had no recorded population, and no votes ever reported. Usually, precincts in rural areas with a tiny population included their results in neighboring precincts for each election, likely because of shared polling places between the precincts. 

To fix this, we assigned zeros to the election results for the non-reporting precincts. For precincts that jointly reported election results, we joined precincts in the shape file if they were neighboring. If precincts reported jointly but were not neighboring, we assigned votes to the precinct closer to more populated areas. 


\vspace{0.75in}
\section{Results for New Hampshire}

\subsection{New Hampshire precinct political geography}
\begin{figure}[!h]
\centerline{\includegraphics[width=\textwidth]{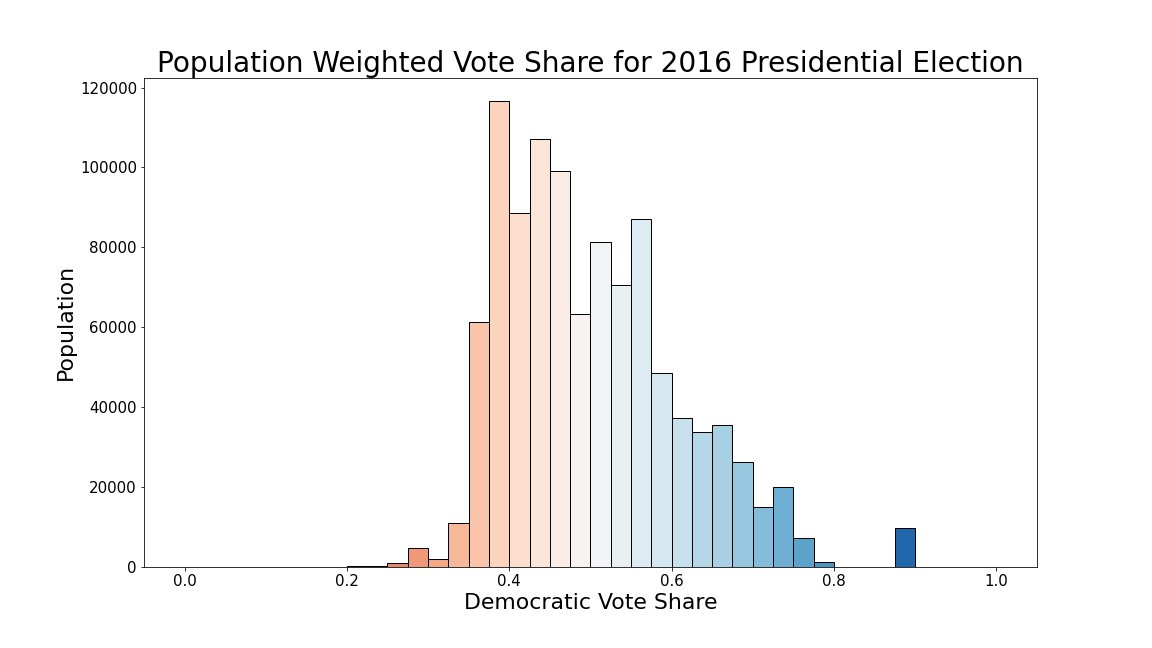}}
\caption{Histogram of Democratic vote share by precinct for the 2016 Presidential election, weighted by population.}
\label{NHPrecinctHist}
\end{figure}

Figure \ref{NHPrecinctHist} shows a population-weighted histogram of the vote share among precincts for the 2016 Presidential election. The height of the bars corresponds to how many people live in a precinct with that Democratic vote share. During this election, New Hampshire was skewed slightly Democratic. More people lived in Republican-favoring precincts, but since there are more highly Democratic precincts ($\geq \frac{2}{3}$) than highly Republican precincts ($\leq \frac{1}{3}$) the state was overall blue. These heavily Democratic precincts create the potential for districts to be packed with Democratic voters more easily than packed with Republican voters, making it somewhat easier for plans that favor Republicans to be created. The dark blue bar on the far right represents the area around Durham, where the University of New Hampshire is.

\subsection{Wins by party}
Using data from the 2016 Presidential Election, we used the GerryChain algorithm to generate 1000 possible congressional redistricting plans. For each plan, we stored whether Democrats would win 0, 1 or 2 of the two possible congressional seats. In 92\% of the plans, one Democrats and one Republican each won a seat (1R 1D). In the remaining 8\% of plans, Democrats won both seats (2D 0R). In the 2018 Election to the US House of Representatives, two Democrats (Chris Pappas and Ann McLane Kuster) won the two districts' seats, with 53.56\% and 55.54\% of the votes, respectively. Making the (admittedly implausible) assumption that voters cast their ballot for the same party as they did in the 2016 Presidential Election, we can conclude that 1R 1D is much more expected than the 2D 0R result that actually occurred; perhaps voters favored the Democratic representative candidates more than the Democratic presidential candidate.

Similarly, we generated 1000 possible congressional redistricting plans using data from the same 2016 Presidential Election, but rather than looking at Democrat (D) versus Republican (R) votes, we looked at Democrat versus Republican combined with Independent (I) votes, and Republican versus Democrat combined with Independent votes. With these results, if all voters who chose an Independent candidate instead chose a Republican, in 1000 possible plans generated by Gerry Chain, 95\% of plans resulted in 0D 2R, 5\% were 1R 1D, and zero plans resulted in 2D. In contrast, if all voters who chose an Independent candidate chose a Democrat instead, in 1000 possible plans, 0 plans resulted in 2R 0D, 2\% of plans resulted in 1R 1D, and 98\% of plans resulted in 2D 0R. This is in line with our explanation in \S\ref{nh-intro} about the importance of Independent voters in New Hampshire elections.

\subsection{Assessing fairness of the enacted plan}
\begin{figure}[!h]
\centerline{\includegraphics[scale=0.4]{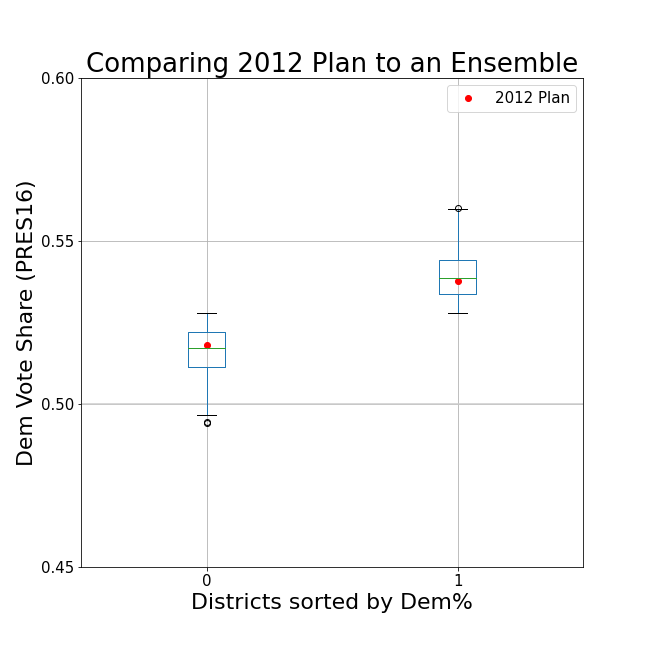}}
\caption{Box and whiskers plot showing how Democratic vote shares from the 2012 enacted plan compared to an ensemble of 1000 plans.}
\label{NHBoxPlot}
\end{figure}

Figure \ref{NHBoxPlot} shows where Democratic vote shares for the 2012 plan fall when compared to an ensemble of 1000 ReCom plans, using the 2016 Presidential election. The box in this plot shows the range of the upper and lower quartile, the center line is the median, and the red dot is the measure for the enacted plan. From this figure, we can see that the 2012 plan does not appear to advantage either party. 

Figure \ref{NHBadEG} highlighted a failure of the efficiency gap to give meaningful information in two-district elections: When a single party wins both of the districts, the efficiency gap is always the same, which leads to a large spike on the histogram near $\frac{1}{2}$ or $-\frac{1}{2}$ that makes it difficult to see the more fine variation in results near 0.   
Figure \ref{NHEG} is essentially a zoomed-in view of the area around 0, this time for an ensemble of 1000 ReCom plans using the 2016 Senate election. For this graph, we restricted results to only show efficiency gaps for plans where each party wins one district. We chose the 2016 Senate election because for the majority of ReCom plans, each party won one district. 

A positive efficiency gap benefits Democrats, so because the center of the histogram is slightly left of 0, it appears that the political geography in New Hampshire slightly benefits Republicans; this is generally true in most of the U.S. The dotted red line represents the efficiency gap of the plan enacted in 2012, and it falls right in the middle of these plans, so the 2012 plan does not appear to be gerrymandered to advantage either party.

\begin{figure}[!h]
\centerline{\includegraphics[scale=0.4]{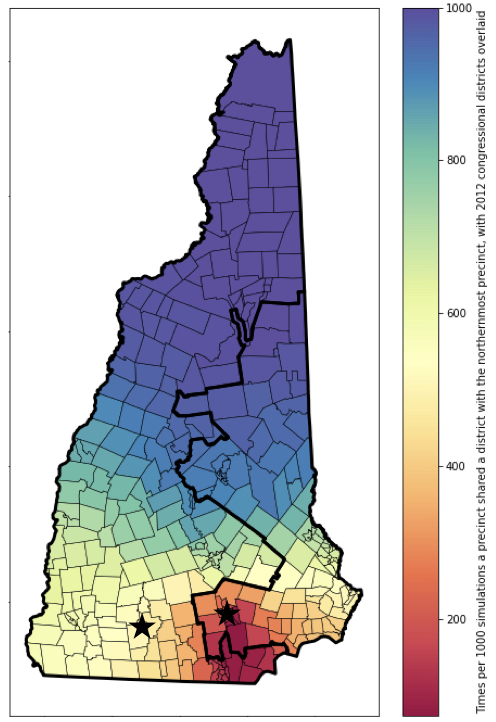}}
\caption{A heat map showing the proportion of plans in which a precinct shared a district with the northernmost precinct (Pittsburg). The 2012 enacted division boundary is outlined, and the hometowns of the incumbent representatives are indicated with stars.}
\label{NHNorthernmost}
\end{figure}

Figure \ref{NHNorthernmost} is a heat map that shows the number of times out of our 1000 ReCom plans a precinct was in the same district as Pittsburg, the precinct comprising the northernmost tip of New Hampshire.\footnote{D.D. thanks Jonathan Mattingly for suggesting this analysis.} The purple region shows precincts that are almost always grouped with Pittsburg; the red region shows precincts that are almost never grouped with Pittsburg; and the yellow region shows precincts that are grouped with Pittsburg about half the time. This means that an ``average'' dividing line between districts is in the yellow region.

Figure \ref{NHNorthernmost} tells us to expect that a districting plan for New Hampshire would (1) keep the entire northern part of the state (purple) together,  (2) keep the urban southern area of Manchester, Nashua and Londonderry (red) together, and (3) divide horizontally between Concord in the center of the state and Manchester just south of it (yellow). In fact, the enacted division does none of these things, and thus is completely contrary to what we expect.

In the figure, the 2012 districting division boundary is shown in black. It runs essentially perpendicular to the places that are most likely to be cut (in yellow) and is much longer than any division line in the yellow region would be. Additionally, both districts include precincts of every color, which is exactly the opposite of what we would expect. 

One possible explanation is that the 2012 plan was chosen to split the largest cities in New Hampshire, which are all located near the red area, into two different districts. Another explanation is \emph{incumbent gerrymandering}: The two incumbent representatives at the time when the map was drawn lived in the precincts indicated by black stars. If the division line went through the yellow region, as we expect, they would both be in the same district and would have to compete; with the division line snaking around their hometowns, they do not have to compete. As it turns out, \emph{both} representatives (Frank Guinta and Charles Bass) were defeated in 2013.

\subsection{County splits}
As explained in \S\ref{gerry-background}, sometimes states have requirements for districting plans, stricter than the equal-population requirement mandated by the U.S. Constitution. A common requirement is to avoid splitting counties and other political units as much as possible. New Hampshire does not have this requirement, but it seems to us that whoever drew the map did take counties into account: Near the top of the division line (see Figures \ref{fig:current_districts} or \ref{NHNorthernmost}) there is a little tendril snaking towards the northwest. This might look like evidence of gerrymandering, but in fact it's just the shape of Carroll county. Thus, we concluded that preserving county boundaries may have been a guiding principle of the map constructor, even though it is not in the law.

\begin{figure}[!h]
\centerline{\includegraphics[scale=1]{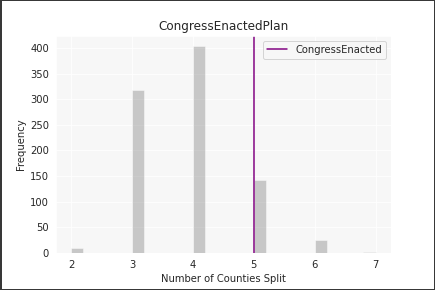}}
\caption{Number of New Hampshire's 10 counties that are split in our ensemble (grey) and in the enacted plan (purple)}
\label{NHSplitCounties}
\end{figure}

Looking at the 2012 congressional plan for New Hampshire, we see that counties are split 6 times, splitting 5 different counties. Considering that New Hampshire only has 10 counties, it is striking that half of the counties have voters in two different congressional districts; it is unclear why anyone would want to split so many. When we drew 1000 districting plans using code that minimized splitting counties, we see that nearly 3/4 of plans have fewer than 5 splits (Figure \ref{NHSplitCounties}). 
Thus, the enacted plan's 5 split counties is slightly more than we would expect, but not extremely so.

\section{Results for Maine}

\subsection{Huge differences in partisanship in different elections}
\begin{figure}[h!]
  \includegraphics[height=0.6\textheight]{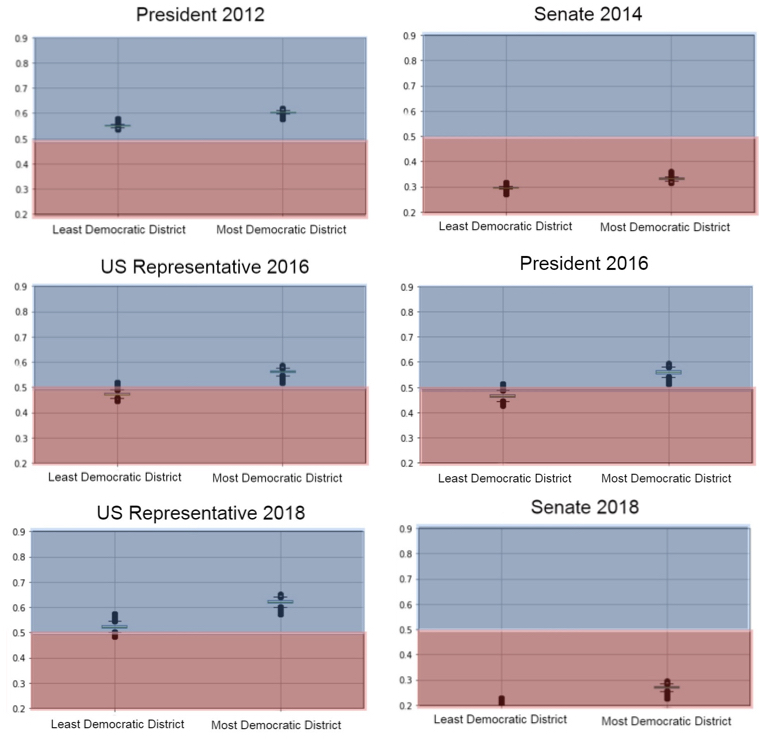}
  \caption{The vertical axis represents the vote share for the Democratic party. Each black dot represents election results for an ensemble of districting plans. Note the very different partisanship between the 2018 Representative and Senate elections (bottom row), which were on the same ballot.}
  \label{fig:voteshare_dem}
\end{figure}

Compared to Maine, New Hampshire's political geography was very easy to analyze. The vote shares earned by the Democrat and Republican parties were reasonably consistent, with each party earning between 40\% and 60\% of the vote in each election and Independent candidates earning roughly 5\% of the vote. This means that in New Hampshire, one could reasonably identify voters or precincts as e.g. ``voting Republican.'' This was not the case in Maine, where (1) on a given ballot the vote shares for different candidates of the same party varied widely (see Figure \ref{fig:voteshare_dem}), (2) independent candidates got a much larger vote share, sometimes enough to win elections, and (3) incumbency effects mean that the identity of the candidate outweighed their party affiliation.

For example, in Figure \ref{fig:voteshare_dem} we can see that while elections sourced from Representative and Presidential elections present a similar picture, with Democratic vote shares at roughly 50\%, elections sourced from Senate elections are heavily skewed  Republican, with Democratic vote shares around only 25\%. In 2014, the heavy Republican skew is due to the incumbency bias for Susan Collins (R), who has been Maine's U.S. Senator since 1996, continuously serving since before all but one of the authors of this paper were born. 

For the Senate 2018 election, the Republican skew is because former Governor Angus King won running as an independent, so plotting only Democrat and Republican votes is misleading since many of the votes cast were for an independent. 
This is a structural limitation with standard measures of gerrymandering, which are defined based on a two-party system and have no clear generalization to a system with candidates from more than two parties.

\subsection{Assessing fairness of the enacted plan}
\begin{figure}[!h]
  \includegraphics[width=0.7\textwidth]{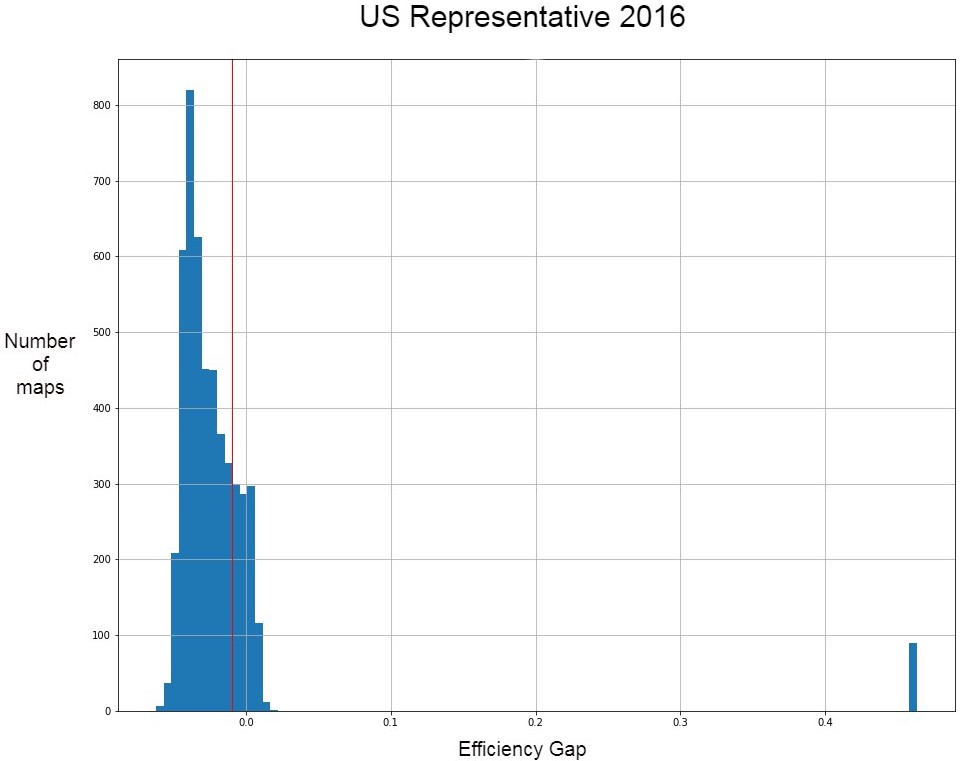}
  \caption{The efficiency gap in Maine's 2016 U.S. Representative election, for an ensemble 5000 legal maps. The red line represents the efficiency gap for Maine's current plan. Because it falls in the middle of the ensemble, this measure does not suggest partisan gerrymandering.}
  \label{fig:rep16_eg}
\end{figure}

As discussed in \S\ref{sec:measures}, the \emph{efficiency gap} measures the extent to which one party ``wastes'' more votes than the other in a given election. As explained above, the efficiency gap has limitations, especially for two-district elections, but by looking at the efficiency gaps across an ensemble of maps, we can compare them to the current plan’s efficiency gap to see how probable or improbable it is that the map was drawn with the intention to waste votes. In Figure \ref{fig:rep16_eg}, we can see that the efficiency gap of the enacted plan (the red line) is similar to the ensemble of 5000 simulated maps. 

\begin{figure}[h!]
  \includegraphics[scale=.4]{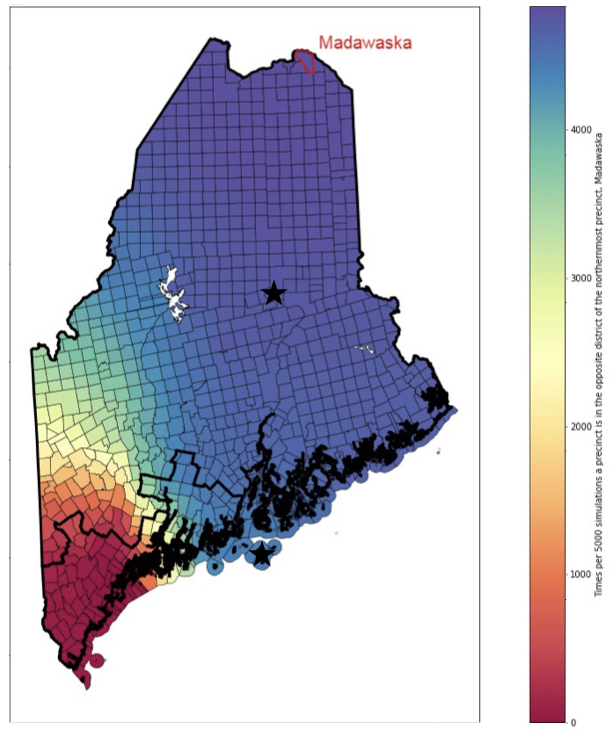}
  \caption{A heat map showing which precincts are most often grouped with Madawaska in a district. The black lines show Maine's current districting plan.}
  \label{fig:maine_heatmap}
\end{figure}

Because of issues using historical election data as a proxy for partisanship in Maine, we analyzed aspects of the plan that do not rely on partisanship. As we did for New Hampshire, we constructed a heat map (Figure \ref{fig:maine_heatmap}) showing how likely each precinct was to be in the same district as Maine's northernmost precinct (Madawaska). In our ensemble of 5000 legal districts, purple precincts were nearly always grouped together, red precincts were nearly always grouped together, and the division line tended to be in the yellow strip.

As in New Hampshire (Figure \ref{NHNorthernmost}), \begin{enumerate}
\item[(1)] both of the enacted districts contain precincts of every color, which our analysis shows to be very unlikely to happen by chance when creating equal-population districts; 
\item[(2)] the enacted dividing line is roughly perpendicular to and much longer than a division line passing through the yellow region, and 
\item[(3)] if the division line did pass through the yellow region, both of the incumbent representatives in 2011, whose hometowns are indicated with black stars, would have ended up in the same region competing against each other, while in the enacted plan they are in different districts. This suggests incumbent gerrymandering.
\end{enumerate}
Unlike in New Hampshire, in this case both of Maine's incumbent representatives (Mike Michaud and Chellie Pingre) were re-elected under this map.

\section*{Acknowledgments}
This work took place during the Districting Data REU in summer 2020, organized by Diana Davis. We thank the Metric Geometry and Gerrymandering Group for their GeoData BootCamp in June 2020, where D.D. and J.L. learned to do this work. D.D., J.L., and D.T. were supported by full or partial summer science research stipends from Swarthmore College.

We thank the other members of the Districting Data REU for their collaboration, technical help, and friendship during this project: 
Ellen Adams, 
Jiahua Chen,
Kevin Choi, 
Kenny Gwon,
William Hoganson, 
Vinay Keefe, 
Keegan McKenna, 
Hulices Murillo, 
Simon Moore, 
Celia Parts, 
Ethan Rothenberg,
Poonam Sahoo, 
Parker Snipes, 
Justin Snyder, 
Theo Uy,
Madeleine Ward,
Jasmine Xie,
Kevin Xin, and
Alicia Yang.

Collecting and cleaning the data for this kind of work is a time-consuming, tedious process that cannot be automated. That is why our team is so large: To complete this project in three months, we needed a lot of person power. A good estimate is that to collect and clean the data required to do this kind of analysis for one state takes between 100 and 1000 person-hours, depending on the size of the state, the condition of the state's data, and the expertise of the people collecting it.

\section*{Interactive supplements}

\subsection*{Draw your own districts} Now that we have collected the data for New Hampshire and Maine, we can use that data not only to analyze past maps (which we did here) and future maps (which we plan to do in 2021 or 2022), but also individual custom maps. Go to \url{https://districtr.org/new} and click on the state of Maine, and then paint your own two-district plan. You can then analyze the simulated results of past elections under your districting plan, and the demographics of your two districts. 

Currently, you can only do this for Maine, but we hope that our data for New Hampshire will be added soon.

\subsection*{Talks} We gave two talks about our work, which you can watch \href{https://www.youtube.com/watch?v=CTR8OGZXIAE}{here} and \href{https://www.youtube.com/watch?v=Ai0z3t2eQ_8}{here}.

Diana Davis \verb+ddavis@exeter.edu+ \\ Phillips Exeter Academy, Department of Mathematics, 20 Main Street, Exeter NH 03833

Sara Asgari \verb+sasgari1@swarthmore.edu+ \\
Quinn Basewitz \verb+qbasewi1@swarthmore.edu+ \\
Ethan Bergmann \verb+ebergma1@swarthmore.edu+ \\
Jack Brogsol \verb+@swarthmore.edu+ \\
Nathaniel Cox \verb+@swarthmore.edu+ \\
Martina Kampel \verb+mkampel1@swarthmore.edu+ \\
Becca Keating \verb+rkeatin1@swarthmore.edu+ \\
Katie Knox \verb+kknox1@swarthmore.edu+ \\
Angus Lam \verb+anguslam852@gmail.com+ \\
Jorge Lopez-Nava \verb+jlopezn1@swarthmore.edu+ \\
Jennifer Paige \verb+jpaige1@swarthmore.edu+ \\
Nathan Pitock \verb+npitock1@swarthmore.edu+ \\
Victoria Song \verb+xsong3@swarthmore.edu+ \\
Dylan Torrance \verb+dtorran1@swarthmore.edu+ \\
Swarthmore College, Department of Mathematics and Statistics\\ 500 College Avenue, Swarthmore PA 19081

\end{document}
